\documentclass[fleqn,reqno,11pt,a4paper,final]{amsart}

\usepackage[a4paper,left=30mm,right=30mm,top=30mm,bottom=30mm,marginpar=20mm]{geometry} 
\usepackage{amsmath}
\usepackage{amssymb}
\usepackage{amsthm}
\usepackage{amscd}
\usepackage[ansinew]{inputenc}
\usepackage{cite}
\usepackage{bbm}
\usepackage{color}
\usepackage[english=american]{csquotes}
\usepackage[final]{graphicx}
\usepackage{hyperref}
\usepackage{calc}
\usepackage{mathptmx}

\graphicspath{{../Pictures/}}

\newtheoremstyle{thmlemcorr}{10pt}{10pt}{\itshape}{}{\bfseries}{.}{10pt}{{\thmname{#1}\thmnumber{ #2}\thmnote{ (#3)}}}
\newtheoremstyle{thmlemcorr*}{10pt}{10pt}{\itshape}{}{\bfseries}{.}\newline{{\thmname{#1}\thmnumber{ #2}\thmnote{ (#3)}}}
\newtheoremstyle{remexample}{10pt}{10pt}{}{}{\bfseries}{.}{10pt}{{\thmname{#1}\thmnumber{ #2}\thmnote{ (#3)}}}
\newtheoremstyle{ass}{10pt}{10pt}{}{}{\bfseries}{.}{10pt}{{\thmname{#1}\thmnumber{ A#2}\thmnote{ (#3)}}}

\theoremstyle{thmlemcorr}
\newtheorem{theorem}{Theorem}

\newtheorem{corollary}[theorem]{Corollary}
\newtheorem{proposition}[theorem]{Proposition}

\theoremstyle{thmlemcorr*}
\newtheorem{theorem*}{Theorem}
\newtheorem{lemma*}[theorem]{Lemma}
\newtheorem{corollary*}[theorem]{Corollary}
\newtheorem{proposition*}[theorem]{Proposition}
\newtheorem{problem*}[theorem]{Problem}
\newtheorem{conjecture*}[theorem]{Conjecture}
\newtheorem{definition*}[theorem]{Definition}

\theoremstyle{remexample}

\theoremstyle{ass}

\def\XXint#1#2#3{{\setbox0=\hbox{$#1{#2#3}{\int}$} 
\vcenter{\hbox{$#2#3$}}\kern-.5\wd0}}

\renewcommand{\epsilon}{\varepsilon}
\renewcommand{\phi}{\varphi}

\begin{document}


\title[Inviscid symmetry breaking with non-increasing energy]{Inviscid symmetry breaking with non-increasing energy}

\author{Emil Wiedemann}
\address{\textit{Emil Wiedemann:} Department of Mathematics, University of British Columbia, and Pacific Institute for the Mathematical Sciences, Vancouver, B.C., Canada V6T 1Z2.}
\email{emil@math.ubc.ca}

\begin{abstract}
In a recent article, C. Bardos et. al. constructed weak solutions of the three-dimensional incompressible Euler equations which emerge from two-dimensional initial data yet become fully three-dimensional at positive times. They asked whether such symmetry-breaking solutions could also be constructed under the additional condition that they should have non-increasing energy. In this note, we give a positive answer to this question and show that such a construction is possible for a large class of initial data. We use convex integration techniques as developed by De Lellis--Sz\'{e}kelyhidi. 
\vspace{4pt}


\vspace{4pt}

\end{abstract}


\hypersetup{
  pdfauthor = {Emil Wiedemann (University of British Columbia and PIMS)},
  pdftitle = {},
  pdfsubject = {},
  pdfkeywords = {}
}


\maketitle




\section{Introduction}
In their recent paper~\cite{bardosetal}, C. Bardos et. al. show, among other results, that there exist weak solutions of the unforced three-dimensional incompressible Euler equations that have two-dimensional initial data, but spontaneously become fully three-dimensional for positive times. They refer to this effect as an instance of ``inviscid symmetry breaking'' and show that there is no ``viscous symmetry breaking'', i.e. a Leray--Hopf solution of the Navier-Stokes equations with two-dimensional initial data remains two-dimensional for all time. Here, by a two-dimensional vectorfield $v:\mathbb{R}^3\to\mathbb{R}^3$ we mean a vectorfield of the form
\begin{equation*}
v(x_1,x_2,x_3)=\left(v_1(x_1,x_2),v_2(x_1,x_2),0\right).
\end{equation*}  
However, the symmetry-breaking solutions of~\cite{bardosetal} are not \emph{dissipative} in the sense of P.-L. Lions~\cite{lions} because their kinetic energy $E(t)=\frac{1}{2}\int|v(t,x)|^2dx$ increases, at least up to some positive time. The authors therefore ask whether there are examples of dissipative weak solutions of the Euler equations which exhibit inviscid symmetry breaking. In fact, such an example is given in~\cite{vortexpaper} for the case of a flat vortex sheet. The aim of this note is to show more generally that the set of initial data for which such dissipative symmetry-breaking occurs is dense in the set of two-dimensional initial data with respect to the topology of $L^2$. Our proof relies on the convex integration approach of De Lellis--Sz\'{e}kelyhidi~\cite{euler2}.
 
We here restrict ourselves to the case of periodic boundary conditions and denote by $\mathbb{T}^d$ the $d$-dimensional torus. By a \emph{weak solution} of the Euler equations on $\mathbb{T}^d$ we mean a solution in the sense of distributions as defined e.g. in~\cite{eulerexistence}, and a \emph{dissipative weak solution} is a weak solution whose energy satisfies $E(t)\leq E(0)$ for all times $t$. It is known~\cite{euler2} that dissipative weak solutions on $\mathbb{T}^d$ are also dissipative solutions in the sense of Lions~\cite{lions} (see Section 4.4 therein), hence the terminology. If $f$ is a measurable function on $\mathbb{T}^3$ such that for almost every $x_1,x_2,a,b\in\mathbb{T}$ we have $f(x_1,x_2,a)=f(x_1,x_2,b)$, then we say that $f$ is essentially independent of $x_3$ (or briefly ei-$x_3$), cf. Definition 4.1 in~\cite{bardosetal}. Finally, $C([0,\infty);H_w(\mathbb{T}^3))$ denotes the space of continuous maps from $[0,\infty)$ into the space $H(\mathbb{T}^3)$ of solenoidal vectorfields equipped with the weak $L^2$ topology (cf. again~\cite{eulerexistence}). Our result can then be formulated as follows:   
\begin{theorem}\label{densesymmetry}
Let $v_0=(v_0^1,v_0^2)\in H(\mathbb{T}^2)$, and let $\epsilon>0$. Then there exist infinitely many dissipative weak solutions $v=v(t,x_1,x_2,x_3)\in C([0,\infty);H_w(\mathbb{T}^3))$ of the incompressible Euler equations such that
\begin{equation*}
||v(t=0)-(v_0,0)||_{L^2(\mathbb{T}^3)}\leq\epsilon,
\end{equation*}
$v(t=0)\cdot e_3=0$, and $v(t=0)$ is ei-$x_3$, but $v$ is not ei-$x_3$ for $t>0$.
\end{theorem}
\section{Proof of Theorem \ref{densesymmetry}}

Since $C^{\infty}(\mathbb{T}^2)$ is dense in $L^2(\mathbb{T}^2)$, we may assume without loss of generality that $v_0$ is smooth. Let $\tilde{v}=\tilde{v}(t,x_1,x_2)$ be the unique solution of the two-dimensional Euler equations in $\mathbb{T}^2$ with initial velocity $v_0$. On defining the $2\times2$--matrix field
\begin{equation*}
\tilde{u}=\tilde{v}\otimes \tilde{v}-\frac{|\tilde{v}|^2}{2}I_2,
\end{equation*}
where $I_2$ denotes the $2\times2$--identity matrix, we obtain a smooth, space-periodic pressure $\tilde{q}$ such that the triplet $(\tilde{v},\tilde{u},\tilde{q})$ satisfies the linear system (1) from Theorem 2 in \cite{eulerexistence}:
\begin{equation}\label{linear}
\begin{aligned}
\partial_t\tilde{v}+\operatorname{div}\tilde{u}+\nabla\tilde{q}&=0\\
\operatorname{div}\tilde{v}&=0.
\end{aligned}
\end{equation}
In~\cite{euler2}, the generalised energy density $e$ is defined by
\begin{equation*}
e(v,u):=\frac{d}{2}\lambda_{max}(v\otimes v-u)
\end{equation*}
for a $d$-dimensional vector $v$ and a symmetric trace-free $d\times d$--matrix $u$, where $\lambda_{max}$ denotes the largest eigenvalue. We have in our situation that
\begin{equation*}
e(\tilde{v}(t,x),\tilde{u}(t,x))=\frac{|\tilde{v}(t,x)|^2}{2}.
\end{equation*}
Now, for some $\eta>0$, choose a function $\bar{e}_0\in C((0,\infty)\times \mathbb{T}^2))\cap C_b([0,\infty);L^1(\mathbb{T}^2))$ with
\begin{equation*}
\frac{|\tilde{v}(t,x)|^2}{2}<\bar{e}_0(t,x)\leq\frac{|\tilde{v}(t,x)|^2}{2}+\eta\hspace{0.3cm}\text{for all $x,t$,}
\end{equation*}
for instance by setting
\begin{equation}\label{defe_0}
\bar{e}_0(t,x)=\frac{|\tilde{v}(t,x)|^2}{2}+\frac{\eta}{1+t}.
\end{equation}

The following proposition follows easily from Proposition 22 in~\cite{euleryoung}, which itself is an adaptation of Proposition 5 in~\cite{euler2} (we omit details):
\begin{proposition}
If in~\eqref{defe_0} $\eta>0$ is chosen sufficiently small, then there exists a triplet $(v',u',q')$ solving~\eqref{linear} in $(0,\infty)\times \mathbb{T}^2$ with the following properties:
\begin{equation*}
(v',u',q')\in C^{\infty}\left((0,\infty)\times \mathbb{T}^2\right),\hspace{0.3cm}v'\in C\left([0,\infty);H_w(\mathbb{T}^2)\right),
\end{equation*}
$u'$ takes values in the space of symmetric trace-free matrices,
\begin{equation*}
e\left(v'(t,x),u'(t,x)\right)<\bar{e}_0(t,x)\hspace{0.3cm}\text{for all $t>0$ and $x\in \mathbb{T}^2$},
\end{equation*}
\begin{equation}\label{prop3}
\frac{1}{2}|v'(0,x)|^2=\bar{e}_0(0,x)\hspace{0.3cm}\text{for a.e. $x\in \mathbb{T}^2$},
\end{equation}
and
\begin{equation}\label{prop4}
||v'(t=0)-\tilde{v}(t=0)||_{L^2(\mathbb{T}^2)}\leq\epsilon.
\end{equation}
\end{proposition}
With this proposition at hand, we define the three-dimensional triplet $(\bar{v},\bar{u},\bar{q})$ by $\bar{v}=(v',0)$,
\begin{equation*}
\bar{u}=\left(\begin{array}{ccc}
u'_{11}+\frac{1}{3}e(v',u')&u'_{12}&0\\
u'_{21}&u'_{22}+\frac{1}{3}e(v',u')&0\\
0&0&-\frac{2}{3}e(v',u')
\end{array}\right)
\end{equation*}
and $\bar{q}=q'+\frac{1}{3}e(v',u')$. An elementary calculation yields that $(\bar{v},\bar{u},\bar{q})$ is a solution of~\eqref{linear} which is smooth for $t>0$ and such that $\bar{v}\in C([0,\infty);H_w(\mathbb{T}^3))$, $\bar{u}$ is symmetric and trace-free, and  
\begin{equation}\label{subsolutionineq}
e(\bar{v}(t,x),\bar{u}(t,x))=e(v'(t,x),u'(t,x))<\bar{e}_0(t,x)\hspace{0.3cm}\text{for all $t>0$ and $x\in \mathbb{T}^3$.}
\end{equation}
(Note that we abused notation in denoting by $e$ both the two-dimensional and the three-dimensional generalised energy density.) Moreover, $\bar{v}(t=0)$ is ei-$x_3$ with
\begin{equation*}
\frac{1}{2}|\bar{v}(0,x)|^2=\bar{e}_0(0,x_1,x_2)\hspace{0.3cm}\text{for a.e. $x\in \mathbb{T}^3$.}
\end{equation*}
Let now 
\begin{equation}\label{defe}
\bar{e}(t,x)=\bar{e}(t,x_1,x_2,x_3)=\bar{e}_0(t,x_1,x_2)+\eta\frac{t}{1+t}\sin^2(2\pi x_3).
\end{equation}
Then, by~\eqref{subsolutionineq},
\begin{equation*}
e(\bar{v}(t,x),\bar{u}(t,x))<\bar{e}(t,x)\hspace{0.3cm}\text{for all $t>0$ and $x\in \mathbb{T}^3$},
\end{equation*} 
so that Theorem 2 in \cite{eulerexistence} yields infinitely many weak solutions $v\in C([0,\infty);H_w(\mathbb{T}^3))$ of the incompressible Euler equations with initial data $\bar{v}(0,\cdot)$ and energy density
\begin{equation}\label{finalenergy}
\frac{|v(t,x)|^2}{2}=\bar{e}(t,x)\hspace{0.3cm}\text{for every $t\in(0,\infty)$ and a.e. $x\in \mathbb{T}^3$}.
\end{equation}
Hence $v$ is ei-$x_3$ at time zero, but not for positive times because its energy density is $x_3$-dependent. Moreover, in view of~\eqref{prop4} we have
\begin{equation*}
||v(t=0)-v_0||_{L^2(\mathbb{T}^3)}=||v'(t=0)-\tilde{v}(t=0)||_{L^2(\mathbb{T}^2)}\leq\epsilon,
\end{equation*}
and finally, owing to~\eqref{finalenergy},~\eqref{defe},~\eqref{defe_0}, and~\eqref{prop3}, for every $t>0$ we have
\begin{equation*}
\begin{aligned}
\int_{\mathbb{T}^3}|v(t,x)|^2dx=&\int_{\mathbb{T}^3}2\bar{e}(t,x)dx\\
\leq&\int_{\mathbb{T}^3}2\bar{e}_0(t,x)dx+\frac{2\eta t}{1+t}\\
&=\int_{\mathbb{T}^2}|\tilde{v}(t,x)|^2dx+2\eta\\
&=\int_{\mathbb{T}^2}2\bar{e}_0(0,x)dx=\int_{\mathbb{T}^3}|v(0,x)|^2dx,
\end{aligned}
\end{equation*}
which completes the proof.\qed

\section{Concluding Remarks}
As a consequence of Theorem~\ref{densesymmetry}, we obtain a plethora of dissipative weak solutions which do not arise as viscosity limits corresponding to the same initial data:
\begin{corollary}
For an $L^2$-dense subset of two-dimensional initial data, there exist dissipative weak solutions of the 3-D incompressible Euler equations which are not a vanishing viscosity limit of Leray-Hopf solutions of Navier-Stokes with the same initial data.
\end{corollary}
To see this, combine Theorem \ref{densesymmetry} with Theorem 4.4 in \cite{bardosetal}. Moreover, the observation that viscosity limits do not break the symmetry while some dissipative weak solutions do can be viewed as further evidence in favour of the viscosity limit as a suitable selection principle for weak solutions of the Euler equations (cf.~\cite{shearflow}).



\providecommand{\bysame}{\leavevmode\hbox to3em{\hrulefill}\thinspace}
\providecommand{\MR}{\relax\ifhmode\unskip\space\fi MR }
\providecommand{\MRhref}[2]{%
  \href{http://www.ams.org/mathscinet-getitem?mr=#1}{#2}
}
\providecommand{\href}[2]{#2}

\end{document}